\documentclass[twoside,10pt,a4paper]{newFNLstyle}
\usepackage{graphics}
\usepackage{cite}
\usepackage{amssymb,amsmath}

\begin{document}

\volnumpagesyear{0}{0}{000--000}{2008}
\dates{received date}{revised date}{accepted date}

\title{COMMENTS ON `REVERSE AUCTION: THE LOWEST UNIQUE POSITIVE INTEGER GAME'}

\authorsone{A. P. FLITNEY}
\affiliationone{School of Physics, \\
University of Melbourne}
\mailingone{Parkville, VIC 3010, Australia. \\
	aflitney@unimelb.edu.au}

\maketitle

\markboth{Comments on Reverse Auction}{A. P. Flitney}

\pagestyle{myheadings}

\begin{abstract}
In Zeng {\em et al.} [Fluct.\ Noise Lett.\ {\bf 7} (2007) L439--L447]
the analysis of the lowest unique positive integer game
is simplified by some reasonable assumptions
that make the problem tractable for arbitrary numbers of players.
However, here we show that the solution obtained
for rational players is not a Nash equilibrium
and that a rational utility maximizer with full computational capability
would arrive at a solution with a superior expected payoff.
An exact solution is presented for the three- and four-player cases
and an approximate solution for an arbitrary number of players. 
\end{abstract}

\section{Introduction}
The lowest unique positive integer game can be briefly described as follows:
Each of $n$ players secretly selects an integer $x$ in the range $[1, n]$
with the player selecting the smallest unique integer
receiving a utility of one,
while the other players score nothing.
If there is no lowest unique integer then all players score zero.

In Zeng {\em et al.}~\cite{zeng07}
the assumption is made that
``...a player is indifferent between two strategies conditioned on the other players' choices
and a player will always pick the lowest number.
Arguably this might be too strong.''
In the following section we show the latter assumption is indeed too strong.

We make use of the following game-theoretic concepts:
A {\em strategy profile} is a set of strategies,
one for each player;
a {\em Nash equilibrium} (NE) is a strategy profile
from which no player can improve their payoff by a unilateral change in strategy;
a {\em Pareto optimal} (PO) strategy profile is one
from which no player can improve their result without someone else being worse off.


\section{Nash equilibrium for the $n=3,4$ player cases}
It is easy to show that the solution of Zeng {\em et al} is not a NE for small $n$.
In the three player case,
if Bob and Charles adopt the strategy $(\frac{1}{2}, \frac{1}{2}, 0)$,
where the $i$th number in the parentheses is the probability of selecting the integer $i$,
Alice can maximize her payoff by selecting the strategy $(0,0,1)$.
In this case Alice wins when ever the other two choose the same integer.
Hence her expected payoff is $\frac{1}{2}$,
double that obtained by selecting the strategy $(\frac{1}{2}, \frac{1}{2}, 0)$.
Bob and Charles win in only $\frac{1}{4}$ of the cases.
This result is an (asymmetric) NE.
Given that $\$_{\text{A}} + \$_{\text{B}} + \$_{\text{C}}  = 1$
the result is also PO:
the sum of the payoffs is maximal so no other strategy profile
can give one player a higher payoff without someone else being worse off.

For $n=4$ there is an analogous solution.
If Bob, Charles and Debra play the strategy $(\frac{1}{2}, \frac{1}{2}, 0, 0)$,
Alice's optimal play is to select `3' with probability one.
Then she wins if the others have all selected `1' or all selected `2'.
The expected payoff to all players is $\frac{1}{4}$
and so the equilibrium is fair to all players.
Again this solution is a NE and is PO
with the maximum possible sum of payoffs (one).
For $n > 4$,
the strategy ``always choose `3' '' is no longer optimal against
a group of players choosing a mixed `1' or `2' strategy
and there is no simple analogue to the above NE strategy profiles.

Asymmetric strategy profiles such as those given above are difficult to realize in practice
since in the absence of communication it is not possible to decide on who plays the odd strategy.
We will now search for a symmetric NE strategy profile
where all the players choose the same (mixed) strategy.
Suppose all players but Alice choose the strategy $(p_1, p_2, \ldots, p_n)$,
while Alice plays $(\pi_1, \pi_2, \ldots, \pi_n)$,
with the normalization conditions\footnote{
In addition, if Alice picks either $n$ or $n-1$
she can only win if all the other players have chosen the same integer.
This will mean that for the NE strategy $\pi_{n-1} = \pi_n$.
However, in the following analysis we shall not make use of this relation.}
$\Sigma p_i = \Sigma \pi_i = 1$.
In the end we will set $\pi_i = p_i \; \forall i$
to give a symmetric strategy profile.
For Alice's strategy to yield her maximum payoff
(given the others' strategies)
it is necessary, though not sufficient, for $d\$_{\text{A}}/d\pi_i = 0, \; \forall i$.

Using the normalization conditions to substitute for $p_n$ and $\pi_n$,
we can write Alice's expected winnings as
\begin{align}
\label{eq:general}
\$_{\text{A}} &= \pi_1 (1-p_1)^{n-1} + \left(1 - \sum_{k=1}^{n-1} \pi_k \right) \, \sum_{j=1}^{n-1} p_j^{n-1} \notag \\
	  & \quad + \: \sum_{i=2}^{n-1} \pi_i
		\left[ (1 - \sum_{j=1}^i p_j)^{n-1} \,+\, \sum_{j=1}^{i-1} p_j^{n-1} \right].
\end{align}
By differentiating with respect to each of the $\pi_i$
and setting the result equal to zero
$n-1$ non-linear coupled equations in the $n-1$ variables $p_1, \ldots p_{n-1}$ are obtained.
Amongst the simultaneous solutions of these equations will be
one that is maximal for Alice.
By setting $\pi_i = p_i \; \forall \, i$
we obtain a strategy that is maximal for all players and is thus a NE.
We note that the derivatives of $\$_{\text{A}}$ do not involve the $\pi_i$.

For the case of $n=3$ we have
\begin{equation}
\$ = \pi_1 (1-p_1)^2 \:+\: \pi_2 \, (p_1^2 + (1-p_1-p_2)^2) \:+\: (1-\pi_1-\pi_2)(p_1^2 + p_2^2),
\end{equation}
(where the subscript $\text{A}$ has been dropped for simplicity) resulting in
\begin{subequations}
\begin{align}
\frac{d\$}{d\pi_1} &= 1 - 2 p_1 - p_2^2 \\
\frac{d\$}{d\pi_2} &= 1 - 2 p_1 +p_1^2 - 2 p_2 + 2 p_1 p_2.
\end{align}
\end{subequations}
This has the unique (for the physical range of $p_1, p_2$) solution 
\begin{equation}
\label{eq:NE3}
p_1 = 2 \sqrt{3} - 3, \qquad
p_2 = 2 - \sqrt{3}.
\end{equation}
We note that $p_2 = 1 - p_1 - p_2$,
as observed in the earlier footnote.
When Bob and Charles play the strategy (\ref{eq:NE3}),
that is, when they select `1' with probability $2 \sqrt{3} - 3 \approx 0.464$
and `2' or `3' each with probability $2 - \sqrt{3} \approx 0.268$,
Alice's payoff is independent of her strategy.
The game being symmetric,
the same is true for any of the players when the other two choose (\ref{eq:NE3}).
Thus, no player can improve their strategy by a unilateral change in strategy,
demonstrating that (\ref{eq:NE3}) is a NE.
When all players select this strategy,
the expected payoff to each is $4(7 - 4 \sqrt{3}) \approx 0.287$,
which is higher than the payoff of $0.25$ that results
when each player selects only between `1' or `2' with equal probability,
the ``rational'' player result of Ref.~\cite{zeng07}.
It is interesting,
and some what anti-intuitive,
that the solution involves a non-zero value for $p_3 = 1 - p_1 - p_2$
since `3' can never be the lowest integer,
though it can be the only unique integer.

Proceeding in the same manner for $n=4$,
(\ref{eq:general}) reduces to
\begin{align}
\$ &= \pi_1 (1-p_1)^3 \:+\: \pi_2 \, [p_1^3 + (1-p_1-p_2)^3] \notag \\
	&\quad +\: \pi_3 \, [p_1^3 + p_2^3 + (1-p_1-p_2-p_3)^3]
	\:+\: (1-\pi_1-\pi_2-\pi_3)(p_1^3 + p_2^3 + p_3^3).
\end{align}
Differentiating with respect to each of $\pi_1, \pi_2$, and $\pi_3$
and setting the results equal to zero
gives the unique (physical) solution 
\begin{equation}
\label{eq:NE4}
p_1 \approx 0.488, \qquad
p_2 \approx 0.250, \qquad
p_3 \approx 0.131,
\end{equation}
again with the relationship $p_3 = 1-p_1-p_2-p_3$.
The exact values for the $p_i$ are complicated and unilluminating.
The payoff to each player when they all choose the strategy (\ref{eq:NE4}),
that is, when each player selects `1' with probability $\approx 0.488$,
`2' with probability $\approx 0.250$ and `3' or `4' each with probability $\approx 0.131$,
is approximately 0.134.
This is higher than that obtainable if all the players simply select between `1' and `2' (0.125).
Again, when three players choose (\ref{eq:NE4})
the payoff to the fourth player is independent of their strategy,
demonstrating that the strategy profile is a NE.
Note the symmetric mixed strategy NE profiles
have lower average payoffs than the asymmetric ones found earlier.

\section{Approximate solution for an arbitrary number of players}
In general, since we have $n-1$ coupled equations of degree $n-1$,
for $n>5$ no analytic solution will be possible,
and for $n=5$ the solution will be problematic.
By inspection of (\ref{eq:NE3}) and (\ref{eq:NE4}) the mixed strategy with
\begin{equation}
\label{eq:approx}
\pi_i = \frac{1}{2^i} \quad \text{for} \; i < n, \qquad
\pi_n = \pi_{n-1} = \frac{1}{2^{n-1}},
\end{equation}
is an approximation to the symmetric NE solutions for $n=3,4$.
Equation (\ref{eq:approx})
is in keeping with our intuition by
giving higher weights to the selection of smaller integers.
The payoff to each player for $n>2$ if all select (\ref{eq:approx}) is
\begin{equation}
\label{eq:pay_approx}
\$ = \sum_{k=1}^{n-1} \left[ \frac{1}{2^k} \sum_{j=1}^{k} \left( \frac{1}{2^j} \right)^{n-1} \right]
		\:+\: \frac{1}{2^{n-1}} \sum_{j=1}^{n-1} \left( \frac{1}{2^j} \right)^{n-1}.
\end{equation}
For $n=3$ the payoff is $\frac{9}{32} \approx 0.281$
and for $n=4$ it is $\approx 0.133$,
both very close to the values for the exact symmetric NE given in the previous section.
The payoff (\ref{eq:pay_approx}) as a function of $n$ is shown in Table~\ref{tab:approx},
along with the payoffs from Ref.~\cite{zeng07} of $1/2^{n-1}$
and the exact solutions for the $n=3$ and $4$ cases.
The payoff given in Ref.~\cite{zeng07}
is slightly smaller than the payoff given by (\ref{eq:pay_approx})
but will asymptote to it as $n$ increases.

\begin{table}[b]
\begin{center}
\begin{tabular}{|l|cccccc|}
\hline
$n$									& 3     & 4     & 5      & 6      & 7      & 8 \\
\hline
Equation (\ref{eq:pay_approx})	& 0.281 & 0.133 & 0.0645 & 0.0317 & 0.0157 & 0.00784 \\
Reference \cite{zeng07}				& 0.25  & 0.125 & 0.0625 & 0.0313 & 0.0156 & 0.00781 \\
Exact										& 0.287 & 0.134 &        &        &        & \\
\hline
\end{tabular}
\end{center}
\caption{\label{tab:approx} The payoff for the approximate symmetric Nash equilibrium solution
of strategy (\ref{eq:approx}) along with the rational player payoffs from Ref.~\cite{zeng07}
and the exact symmetric Nash equilibrium payoffs
of strategies (\ref{eq:NE3}) and (\ref{eq:NE4}),
for the three- and four-player cases, respectively.
Exact solutions for the other cases have not be calculated.
Payoffs have been rounded to three significant figures.}
\end{table}

\section{Conclusion}
We have found both asymmetric and symmetric NE strategy profiles
for a three- and four-player
lowest unique positive integer game
with payoffs superior to that resulting
from the simplifying assumption of Ref.~\cite{zeng07}.
In particular the assumption that a player will always choose the lowest integer
in a situation where they have a choice
results in a strategy that is
not a NE.
The asymmetric NE are also PO,
and in the case of $n=4$,
is fair to all players.
The symmetric solutions are unique amongst symmetric strategy profiles
but yield a lower payoff than the asymmetric solutions.
Anti-intuitively,
the NE strategy profiles includes a non-zero probability
for selecting the largest integer
since this may be the only unique integer.

For arbitrary $n$,
we propose a simple symmetric strategy profile
with geometrically decreasing probabilities of selecting higher integers.
This gives very close to the payoffs of the exact symmetric NE solutions
for the two cases for which exact solutions were obtained.
The rational player solution of Ref.~\cite{zeng07} is simpler than ours
but gives payoffs slightly smaller.

\section*{Acknowledgment}
Funding was provided by the Australian Research Council
grant number DP0559273.

\end{document}